\documentclass[letterpaper, 10pt, conference]{ieeeconf}
\IEEEoverridecommandlockouts
\usepackage{latexsym, amssymb,amsmath, amsbsy,amsopn, amstext}
\usepackage{graphicx,epsfig,tikz,pgf,hyperref,cancel,color,cite,float,ifpdf}

\def\G{\mathcal{G}}
\def\E{\mathcal{E}}
\def\V{\mathcal{V}}
\def\A{\mathcal{A}}
\def\D{\mathcal{D}}
\def\d{\delta}

\newcommand{\F}{{\mathbb F}_p}

\newtheorem{thm}{Theorem}[section]
\newtheorem{lem}[thm]{Lemma}

\newtheorem{cor}[thm]{Corollary}
\newtheorem{dfn}[thm]{Definition}

\newtheorem{rem}[thm]{Remark}

\begin{document}

\title{\LARGE \bf Leader-following Consensus of Multi-agent Systems over Finite Fields}

\author{Xiangru Xu,\hskip 2mm  Yiguang Hong
\thanks{Xiangru Xu and Yiguang Hong are with Key Laboratory  of Systems and
Control,  Academy of Mathematics and Systems Science, Chinese
Academy of Sciences, Beijing 100190, China. E-mail: {\tt\small
xuxiangru@amss.ac.cn;  yghong@iss.ac.cn.}}
}

\maketitle

\begin{abstract}
The leader-following consensus problem of multi-agent systems over
finite fields $\F$ is considered in this paper. Dynamics of each
agent is governed by a linear equation over $\F$, where a
distributed control protocol is utilized by the followers.
Sufficient and/or necessary conditions on system matrices and graph
weights in $\F$ are provided for the followers  to track the
leader.
\end{abstract}

\section{Introduction}

The distributed control of multi-agent systems have attracted
intensive attentions these years. Various approaches are proposed to
handle different problems for agents with different communication
and dynamic constraints. In most existing literature,
the states of agents and the information exchange between agents are
defined as real numbers or quantized values
\cite{ber03},\cite{kas07}, \cite{lyn96},\cite{sab07}. Recently, a finite field
formalism was proposed to investigate multi-agent systems where the
states of each agent are considered elements of a finite field
\cite{pas13},\cite{sun13}. The states of each agent are updated
iteratively as a weighted sum of the states of its neighbors, where
the operations are performed as modular arithmetic in that field.
Such a system is not only interesting theoretically but also has
advantages such as smaller convergence time and resilience to
communication noises, with applications to quantized control and distributed 
estimation.

Dynamical systems that take values from finite sets are ubiquitous.
Consensus or synchronization of such systems were also widely
investigated, such as quantized consensus, logical consensus,
synchronization of finite automata \cite{fag13},\cite{jun12},
\cite{kas07},\cite{xu13b}. In fact, finite
fields provide one convenient approach to model some of these
systems. In the communication and circuits areas, linear systems
over finite fields have long been studied\cite{els59},
\cite{koe03},\cite{lid96}. In the control community,
Kalman et al developed an algebraic theory for linear systems over
an arbitrary field in the 1960s, by merging automata theory and
module theory \cite{kal69}. However, to the best of our knowledge,
there are few results on the consensus/synchronization of linear
multi-agent systems over finite fields, especially in a distributed
manner. It is very interesting to consider how to achieve desired
collective behaviors through local information for such systems. In
\cite{sun13}, a first-principle approach was proposed to establish a
graph-theoretic characterization of the controllability and
observability problems for linear systems over finite fields. These
results were applied to state placement and information dissemination
of agents whose states are quantized values. In \cite{pas13}, some
sufficient and necessary conditions on network weights and topology
were given for the consensus of a group of agents on finite fields.
It was shown that analyzing tools for real valued multi-agent systems
cannot be applied straightforwardly to these systems in finite fields.

The objective of this paper is to study the leader-following
consensus problem of multi-agent systems over finite fields.
Dynamics of the leader and the followers are governed by linear
equations in a given finite field. For the leader, the equation is
autonomous; for each follower, it has local information input that is a weighted
sum of relative states between itself and its neighbors, where the
operations are done as modular arithmetic. We first formulate the
leader-following consensus problem on finite fields. Then under some
assumptions, we provide sufficient and/or necessary consensus
conditions on system matrices and graph weights. Compared with
existing results on multi-agent systems over finite fields
\cite{pas13},\cite{sun13}, agents considered here have higher order
dynamics and the interaction graphs are directed acyclic and could
possibly be time-varying.

The rest of the paper is organized as follows. In section 2, some
preliminaries on finite fields and linear systems over finite fields
are given. After the problem is formulated in section 3, our main
results are provided in section 4, along with an illustrative
example. In section 5 some conclusions are presented finally.

\section{Preliminaries}

In this section, preliminary knowledge on finite
field, linear system over finite fields and graph theory will be presented for convenience.

\subsection{Finite Field}

\begin{dfn}\label{dfField}
A field is a commutative division ring. Formally, a field $\mathbb{F}$ is a set of elements with addition ($+$) and multiplication ($\cdot$) operations such that the following axioms hold:
\begin{itemize}
  \item Closure under addition and multiplication.\\
  $\forall a,b\in \mathbb{F}$, $a+b\in \mathbb{F}$, $a\cdot b\in \mathbb{F}$.
  \item Associativity of addition and multiplication.\\
  $\forall a,b,c\in \mathbb{F}$, $a+(b+c)=(a+b)+c$, $a\cdot(b\cdot c)=(a\cdot b)\cdot c$.
  \item Commutativity of addition and multiplication.\\
  $\forall a,b\in \mathbb{F}$, $a+b=b+a$, $a\cdot b=b\cdot a$.
  \item Existence of additive and multiplicative identity elements.\\
 $\exists b,c \in \mathbb{F}$, $\forall a\in \mathbb{F}$, $a+b=a$, $a\cdot c=a$.
  \item Existence of additive and multiplicative inverse elements.\\
  $\forall a\in \mathbb{F}$, $\exists b\in \mathbb{F}$ such that $a+b=0$; $\forall a\in \mathbb{F}, a\neq 0$, $\exists c\in \mathbb{F}$ such that $a\cdot c=1$.
  \item Distributivity of multiplication over addition.\\
  $\forall a,b,c\in \mathbb{F}$, $a\cdot(b+c)=(a\cdot b)+(a\cdot c)$.
\end{itemize}
\end{dfn}

The number of elements (or the \emph{order}) of a finite field is $p^n$, where $p$ is a prime number and $n$ is a positive
integer. Therefore, a finite field is denoted as $\mathbb{F}_{p^n}$.
If $n=1$, $F_p\cong \mathbb{Z}/p\mathbb{Z}=\{0,1,...,p-1\}$, and the
addition and multiplication are done by the module $p$ arithmetic.
If $n>1$, $\mathbb{F}_{p^n}\cong \mathbb{F}_p[x]/(f(x))$, where
$f(x)$ is an irreducible polynomial in $\mathbb{F}_p[x]$. In this
study, we consider finite fields $\mathbb{F}_p$ where $p$ is a
prime number.


The finite field $\F$ is not algebraically closed, which means that
not every polynomial with coefficients in $\F$ has a root in $\F$.
Therefore, not all $N\times N$ matrices  have $N$ eigenvalues in
$\F$. This fact makes many eigenvalue-based results such as the PBH
test for controllability (observability) and the consensus
conditions in the real-valued dynamical systems fail
\cite{pas13,sun13}.

In what follows, we denote zero matrix of dimension $m_1\times m_2$  in $\F$ by $\textbf{0}_{m_1\times m_2}$ and omit  the foot indices  if  clear from context; denote zero vector of dimension $m$ in $\F$ by $0_m$.

\subsection{Linear System over Finite Field}

Consider an autonomous dynamical system over $\F$ as follows:
\begin{equation}\label{eqnAutFin}
x(k+1)=Ax(k)
\end{equation}
where $A\in\F^{n\times n}$, $x(k)\in \F^{n\times 1}$.

Suppose that $P_\lambda(A)=det(\lambda I_n-A)$ is the characteristic
polynomial of $A$, which can be factorized in $\F$ as
$P_\lambda(A)=\lambda^sQ(\lambda)$  with $Q(0)\neq 0$. Dynamics of
the system is completely determined by $P_\lambda(A)$.

\begin{lem}\label{lemDynFin}\cite{tol05}
Dynamics of (\ref{eqnAutFin}) is the product of a tree, which
corresponds to the nilpotent part $\lambda^s$, and the cycles, which
correspond to the bijective part $Q(\lambda)$.
\end{lem}

Consider a dynamical system with control over $\F$ as follows:
\begin{equation}\label{eqnLin}
x(k+1)=Ax(k)+Bu(k)
\end{equation}
where $A\in\F^{n\times n}$ and $B\in\F^{n\times m}$.

The controllability indices $c_i(i=1,...,m)$ of $(A,B)$ can be
defined exactly the same way as for real-valued systems \cite{reg04}. If
$\sum_{i=1}^mc_i=n$, then the matrix
$[B,AB,...,A^{n-1}B]\in\F^{n\times nm}$ has (full) rank $n$, and the
system $(A,B)$ is controllable. If $\bar n:=\sum_{i=1}^mc_i<n$, $(A,B)$ is
uncontrollable and can be partitioned into controllable and
uncontrollable parts.

\begin{lem}\cite{reg04}\label{thmUncon}
Consider control system (\ref{eqnLin}) over $\F$. There is a
transformation of state coordinate $x^c=Qx$ with nonsingular matrix
$Q$ such that (\ref{eqnLin}) can be transformed into a system of the
following form:
\begin{equation}\label{eqnCCF2}
x^c(k+1)=\left(
           \begin{array}{cc}
             A^c & A^{cc} \\
             \textbf{0} & A^{uc} \\
           \end{array}
         \right)
x^c(k)+\left(
          \begin{array}{c}
            B^c \\
            \textbf{0} \\
          \end{array}
        \right)
u(k)
\end{equation}
where $\left(
           \begin{array}{cc}
             A^c & A^{cc} \\
             \textbf{0} & A^{uc} \\
           \end{array}
         \right)=QAQ^{-1}$, $\left(
          \begin{array}{c}
            B^c \\
            \textbf{0} \\
          \end{array}
        \right)=QB$, $A^c\in\F^{\bar n\times\bar n},B^c\in\F^{\bar n\times m}$,$(A^c, B^c)$ is controllable and in the control companion form.
%
\end{lem}

\begin{dfn}\label{dfnnil}
A \emph{nilpotent} matrix over $\F$ is a square matrix $A\in
\F^{m\times m}$ such that $A^k=\textbf{0}$ for a
positive integer $k$. The smallest $k$ to satisfy $A^k=\textbf{0}$  is called
the nilpotent \emph{degree} of $A$.
\end{dfn}


\begin{dfn}\label{dfnSta}
The system $(A,B)$ is called \emph{stabilizable} if the uncontrollable subsystem matrix $A^{uc}$ in (\ref{eqnCCF2}) is nilpotent.
\end{dfn}

By Lemma \ref{thmUncon}, it is not hard to see that $(A,B)$ is stabilizable if and only if there is a matrix $K\in\F^{m\times n}$ such that $A+BK$ is nilpotent.

\subsection{Graph Theory}

The information exchange between agents is described by a graph
$\G=\{\V,\E\}$, where $\V=\{1,...,N\}$ is the set of vertices to
represent $N$  agents and $\E\subset\V\times\V$ is the set of edges
to represent the information exchange between agents. If
$(i,j)\in\E$, then agent $j$ can receive information from agent $i$.
The set of \emph{neighbors} of the $i$-th agent is denoted by
$N_i=\{j\in \V|(j,i)\in\E\}$.  In this study the graph considered is
directed, that is, $(i,j)\in\E$ not necessarily implies
$(j,i)\in\E$. If there exists a sequence of nodes $i_1,i_2,...,i_t$
such that $(i_j,i_{j+1})\in\E$ for $j=1,...,t-1$, then the sequence
is called a \emph{path} from node $i_1$ to $i_t$ and the node $i_t$
is called \emph{reachable} from $i_1$. If $i_t=i_1$, then the path
is called a cycle. The union of a set of graphs
$\{\G_1=\{\V_1,\E_1\},...,\G_m=\{\V_m,\E_m\}\}$ is a directed graph
with nodes given by $\cup_{i=1}^m \V_i$ and edge set given by
$\cup_{i=1}^m\E_i$.

Given a finite field $\F$, the \emph{weighted adjacency matrix} of
$\G$ is denoted as $\A=(a_{ij})\in \F^{N\times N}$, where
$a_{ij}=0$ if $(j,i)\notin\E$. Here, ``0''
is the additive identity of $\F$. The \emph{in-degree} of
node $i$ is defined as $d_i=\sum_{j=1}^Na_{ij}$ and the Laplacian
matrix of $\G$ is defined as $L=\D-\A$ where $\D=diag(d_1,...,d_N)$
is the \emph{degree matrix}. A directed graph without cycles is
called a \emph{directed acyclic graph} (DAG). Suppose that $\A$ is
the weighted adjacency matrix of a directed graph $\G$. Then $\G$ is
DAG if and only if there is a permutation matrix $P$ such that $P\A
P^{-1}$ is strictly upper triangular \cite{nic75}.



\section{Problem Statement}
Given a finite field $\F$, let us consider a multi-agent system
consisting of one leader represented by $0$ and $N$ followers
represented by $\{1,...,N\}$. The state of agent $i$ is described by
a column vector of dimension $n$: $x_i=(x_i^1,...,x_i^n)^T$ with
$x_i^s\in\F (i=0,...,N,s=1,...,n)$. The interaction graph describing
the information exchange among the $N+1$ agents is denoted by
$\G=(\V,\E)$, while the subgraph induced by the $N$ followers is
denoted by $\bar\G$. The weighted adjacency matrix
and degree matrix of the $N+1$ agent system are denoted by
$\A=(a_{ij})\in\F^{(N+1)\times (N+1)}$, $\D\in\F^{(N+1)\times
(N+1)}$, respectively. Correspondingly, the induced adjacency submatrix and degree submatrix corresponding to $\bar\G$ are denoted by
$\bar\A\in\F^{N\times N}$ and $\bar\D\in\F^{N\times N}$,
respectively.

Dynamics of the leader is described by a linear equation over $\F$ as follows:
\begin{equation}\label{eqnLeader}
  x_0(t+1) = Ax_0(t)
\end{equation}
where $A\in\F^{n\times n}$.

Dynamics of the $i$-th follower is described by a linear control system over $\F$ as follows:
\begin{equation}\label{eqnFollow}
  x_i(t+1) = Ax_i(t)+bu_i(t)
\end{equation}
where $b=(b_1,...,b_n)^T\in\F^{n\times 1}$ is a column vector and
$u_i(t)\in \F$ is the input. Note that the addition and
multiplication in (\ref{eqnLeader}) and (\ref{eqnFollow}) are
modular arithmetic in $\F$.

In \cite{pas13}, consensus of agents over $\F$ was studied where the state
of each agent is represented by a scalar in $\F$. Consensus is said to be
achieved if all the agents eventually have the same value. Dynamics
of the overall agent network can be described by the autonomous
equation (\ref{eqnLeader}). Unlike the real-valued discrete-time
consensus problem, it was shown in \cite{pas13} that $\G$ is
strongly connected and matrix $A$ is row-stochastic can not
guarantee consensus in $\F$, and tools for analyzing real valued multi-agent systems cannot be applied straightforwardly to these systems in $\F$. Necessary and sufficient conditions were provided in \cite{pas13} for consensus and average consensus, which were much more restrictive compared with the corresponding results for real-valued systems.

In this paper, we consider high-order agents (\ref{eqnLeader}) and (\ref{eqnFollow}) rather than the scalar
agents discussed in \cite{pas13}. Correspondingly, the consensus
problem of (\ref{eqnLeader}) and (\ref{eqnFollow}) is defined as
follows.

\begin{dfn}\label{dfnCon}
The followers (\ref{eqnFollow}) achieve (finite-time) consensus with
the leader (\ref{eqnLeader}) in $\F$ if for any initial state $x_0(0)$,
$x_i(0),i=1,...,N,$ there exists $T\in \mathbb{Z}_+$ such that for
any $s=1,...,n$ and $k\geq T$,
\begin{equation}\label{finiteCon}
x_i^s(k)=x_0^s(k)
\end{equation}
\end{dfn}

\begin{rem}
Actually, it is not hard to show that the \emph{finite-time consensus}
defined in (\ref{finiteCon}) is equivalent to the \emph{asymptotic
consensus} defined as $\lim_{k\rightarrow \infty}x_i^s(k)=x_0^s(k)$.
Therefore, we only say that (\ref{eqnFollow}) achieve \emph{consensus} with (\ref{eqnLeader}). 
\end{rem}


Suppose that the input in (\ref{eqnFollow}) has the following form:
\begin{equation}\label{eqnControl}
  u_i(k) = K\sum_{j=0}^Na_{ij}(x_j(k)-x_i(k))
\end{equation}
where $K\in\F^{1\times n}$ is a constant matrix.

Consensus problem of
real-valued discrete-time multi-agent systems with control
(\ref{eqnControl}) has been intensively investigated for both leadless
and leader-following cases
\cite{hong06},\cite{kin06},\cite{ma10},\cite{mov13},\cite{ni10},\cite{sun07}. However,
existing analyzing tools do not apply straightforwardly to our problem.

In this study, we aim to find sufficient and/or necessary conditions
on system matrices $(A,b)$ and weighted digraph $\G$ to make (\ref{eqnFollow})
achieve consensus with (\ref{eqnLeader}) under control protocol of the form
(\ref{eqnControl}).

\section{Consensus Conditions}

In this section, we give our main results on consensus conditions.

Consider equations (\ref{eqnLeader}) and (\ref{eqnFollow}). Let
$\d_i(k)=x_i(k)-x_0(k)\in\F^{n\times 1} (i=1,...,N)$ and $\d_0(k)=0_n$ for
simplicity. Then consensus condition (\ref{finiteCon}) is equivalent
to the existence of $T\in \mathbb{Z}_+$ such that, for
any $i=1,...,N$ and $k\geq T$,
\begin{equation}\label{connil}
\d_i(k)= 0_n
\end{equation}
For any $j,l\in\{1,...,N\}$, $\d_j(k)-\d_l(k)=x_j(k)-x_l(k)$. Then for any $i=1,...,N$,
\begin{eqnarray*}
  \d_i(k+1)&=& A\d_i(k)+bK\sum_{j=0}^Na_{ij}(x_j(k)-x_i(k)) \\
   &=&A\d_i(k)+bK\sum_{j=0}^Na_{ij}(\d_j(k)-\d_i(k))  \\
   &=&A\d_i(k)+bK[\sum_{j=1}^Na_{ij}\d_j(k)-d_i\d_i(k)]
\end{eqnarray*}

Denote $\d(k)=[\d_1(k)^T,...,\d_N(k)^T]^T$. Then
\begin{equation}\label{eqnErr}
\d(k+1)=[I_N\otimes A+(\bar \A-\bar \D)\otimes bK]\d(k)
\end{equation}

Clearly, condition (\ref{connil}) is equivalent to that $0_{nN}$ is the only
equilibrium of (\ref{eqnErr}). In other words, the matrix
$I_N\otimes A+(\bar \A-\bar \D)\otimes bK$ is nilpotent in $\F$.

In what follows, we assume that the induced subgraph $\bar \G$ is a
directed acyclic graph (DAG), which was actually used in many
existing studies of multi-agent consensus \cite{das02},\cite{shi09},\cite{tan04}. Note also
that DAG is different from the graph topology discussed in
\cite{sun13}, where the topology was assumed to be a spanning tree
(or forest) with self-loops.

If $\bar \G$ is DAG, then there exists a permutation matrix $P$ such
that $P\bar \A P^{-1}$ is a strictly upper-triangular matrix,
denoted as $\hat \A$. Let $\bar D=diag(d_1,...,d_N)$ and $\hat
D=P\bar \D P^{-1}=diag(\hat d_1,...,\hat d_N)$. Then,
$$
\begin{array}{rl}
&(P\otimes I_n)[I_N\otimes A+(\bar \A-\bar \D)\otimes bK](P\otimes I_n)^{-1}\\
=&I_N\otimes A+(\hat \A-\hat \D)\otimes bK
\end{array}
$$

Since $\hat \A$ is strictly upper-triangular and $\hat\D$ is diagonal, $I_N\otimes A+(\hat \A-\hat \D)\otimes bK$ has the following form:
\begin{equation}\label{eqnDIA}
\left(
  \begin{array}{cccc}
    A-\hat d_1 bK & * & * & * \\
    \textbf{0} & A-\hat d_2 bK & * & * \\
    \vdots & \vdots & \ddots & * \\
    \textbf{0} & \cdots & \textbf{0} & A-\hat d_N bK \\
  \end{array}
\right)
\end{equation}


\begin{lem}
If $\bar \G$ is DAG, then condition (\ref{connil}) holds if and only if $A-\hat
d_i bK$ (or equivalently, $A-d_i bK$) is nilpotent in $\F$ for
$i=1,...,N$.
\end{lem}

{\bf Proof:} Note that matrix $M\in\F^{\ell\times \ell}$ is
nilpotent if and only if its characteristic polynomial satisfies
$det(\lambda I-M)=\lambda^\ell$. Because of the upper-triangular block form of \ref{eqnDIA}, $det(\lambda
I_{Nn}-(I_N\otimes A+(\bar \A-\bar \D)\otimes BK))=\Pi_{i=1}^N
det(\lambda I_n-(A-\hat d_i bK))$. Then $det(\lambda
I_{Nn}-(I_N\otimes A+(\bar \A-\bar \D)\otimes BK))=\lambda^{Nn}$ if
and only if $det(\lambda I_n-(A-\hat d_i bK))=\lambda^n$, which is equivalent to that $A-\hat
d_i bK$ (or $A-d_i bK$) is nilpotent for $i=1,...,N$. \hfill$\Box$

Apply Lemma \ref{thmUncon} to (\ref{eqnFollow}). Then there exists an invertible matrix $Q\in\F^{n\times n}$ such that
$$
QAQ^{-1}=A_Q=\left(
           \begin{array}{cc}
             A^c & A^{cc} \\
             \textbf{0} & A^{uc} \\
           \end{array}
         \right),\;
         Qb=b_Q=\left(
           \begin{array}{c}
             b^c \\
             \textbf{0} \\
           \end{array}
         \right)
         $$
where
\begin{equation*}\label{compa}
A^c=\left(
                 \begin{array}{cccccc}
                   0 & 1 & 0 & \cdots & 0 \\
                   0 & 0 & 1 & \cdots & 0 \\
                   \vdots  & \vdots & \vdots & \ddots & \vdots \\
                  0 & 0 & 0  & \cdots & 1 \\
                   a_1 & a_2 & a_3&\cdots &a_{s}\\
                 \end{array}
               \right),
               \;b^c=\left(
               \begin{array}{c}
               0 \\
               0 \\
               \vdots \\
               0 \\
               1 \\
               \end{array}
               \right)
\end{equation*}
and $s$ is the controllability index. Letting $K^c=KQ^{-1}$, we have
equation (\ref{eqnCan}), shown on the next page.

\begin{figure*}[!hbt]
\begin{equation}\label{eqnCan}
\begin{array}{ll}
(I_N\otimes Q)(P\otimes I_n)[I_N\otimes A+(\bar \A-\bar \D)\otimes bK](P\otimes I_n)^{-1}(I_N\otimes Q)^{-1}&\\
=I_N\otimes A_Q+(\hat \A-\hat \D)\otimes b_QK^c&\\
=\left(
  \begin{array}{ccccccc}
      A^c-\hat d_1 b^cK^c & *& * & * & *& * & * \\
      \textbf{0} & A^{uc} & * & * & * & * & * \\
    \textbf{0} & \textbf{0} &A^c-\hat d_2 b^cK^c & *& * & * & *\\
    \textbf{0} & \textbf{0} &\textbf{0} & A^{uc}& * & * & *\\
    \vdots & \vdots & \vdots & \vdots &\ddots & \vdots & \vdots\\
    \textbf{0} & \textbf{0}& \textbf{0}& \cdots&  \textbf{0}&
      A^c-\hat d_N b^cK^c & * \\
      \textbf{0} &\textbf{0} &\textbf{0} &\textbf{0} &\textbf{0} &\textbf{0} &A^{uc}
\end{array}
\right)
\end{array}
\end{equation}
\end{figure*}



For simplicity, we assume that the matrix $A$ is not nilpotent. In
fact, if $A$ is nilpotent, then consensus can be easily achieved by
just letting $K=0$ regardless of $b$ and $\G$.

\subsection{Lemmas On Nilpotent Matrices}
Two lemmas on nilpotent matrices are presented in this subsection for later use.

\begin{lem}\label{lemNil1}
Suppose that matrix $A\in \F^{(n_1+n_2)\times (n_1+n_2)}$ has the
following form
$$
A=\left(
    \begin{array}{cc}
      A_1 & A_2 \\
      \textbf{0}_{n_2\times n_1} & A_3 \\
    \end{array}
  \right)
$$
where $A_1\in \F^{n_1\times n_1},A_3\in \F^{n_2\times n_2}$ are two
nilpotent matrices with nilpotent degrees $k_1$ and $k_2$,
respectively.
Then $A$ is also nilpotent with nilpotent degree upper bounded by $k_1+k_2$. 
\end{lem}

{\bf Proof:}
$$
\begin{array}{rl}
A^{k_1+k_2}&=A^{k_1}A^{k_2}\\
&=
    \left(
      \begin{array}{cc}
        A_1^{k_1} & * \\
        \textbf{0}_{n_2\times n_1} & A_3^{k_1} \\
      \end{array}
    \right)
    \left(
      \begin{array}{cc}
        A_1^{k_2} & * \\
        \textbf{0}_{n_2\times n_1} & A_3^{k_2} \\
      \end{array}
    \right)
\end{array}
$$
If $k_1>k_2$, then
$$
\begin{array}{rl}
A^{k_1+k_2}
    &=
    \left(
      \begin{array}{cc}
        \textbf{0}_{n_1\times n_1} & * \\
        \textbf{0}_{n_2\times n_1} & \textbf{0}_{n_2\times n_2}\\
      \end{array}
    \right)
    \left(
      \begin{array}{cc}
       A_1^{k_2}  & * \\
        \textbf{0}_{n_2\times n_1} & \textbf{0}_{n_2\times n_2}\\
      \end{array}
    \right)\\
    &=\textbf{0}_{(n_1+n_2)\times (n_1+n_2)}
\end{array}
$$
If $k_2\geq k_1$, then
$$
\begin{array}{rl}
A^{k_2+k_1}
    &=
      \left(
      \begin{array}{cc}
        \textbf{0}_{n_1\times n_1} & * \\
        \textbf{0}_{n_2\times n_1} & A_3^{k_1}\\
      \end{array}
    \right)
    \left(
      \begin{array}{cc}
        \textbf{0}_{n_1\times n_1} & * \\
        \textbf{0}_{n_2\times n_1} & \textbf{0}_{n_2\times n_2}\\
      \end{array}
    \right)\\
    &=\textbf{0}_{(n_1+n_2)\times (n_1+n_2)}
\end{array}
$$
The conclusion follows immediately.\hfill$\Box$

Consider a finite set of matrices $\{A^1,...,A^q\}$ over $\F$. Here $A^i=(A^i_{jk})\in \F^{ns\times ns}$ is a block matrix of the following form
\begin{equation}\label{eqnBloNil}
A^i=\left(
  \begin{array}{cccc}
    A_{11}^i & A_{12}^i  & \cdots  & A_{1s}^i  \\
    \textbf{0} & A_{22}^i  & \cdots & A_{2s}^i \\
    \vdots & \vdots & \ddots & \vdots \\
    \textbf{0} & \cdots & \textbf{0} & A_{ss}^i \\
  \end{array}
\right)
\end{equation}
where $A_{jk}^i\in \F^{n\times n}$ and $A_{jj}^{i}=A_{j}$ for some matrices $A_j$, $i=1,...,q$, $j,k=1,...,s$.

\begin{lem}\label{lemNil2}
Consider a finite set of matrices $\{A^1,...,A^q\}$ over $\F$ as shown in
(\ref{eqnBloNil}). Suppose that $A_{j} (j=1,...,s)$ are nilpotent
matrices with respective nilpotent degree $k_j$. Then there exists
an integer $T=\sum_{j=1}^s\tau_j>0$ such that, for any $t\geq T$ and
any sequence $A^{i_1},A^{i_2},...$, we have
\begin{equation*}
A^{i_1}A^{i_2}...A^{i_t}=\textbf{0}_{ns\times ns}
\end{equation*}
where $\tau_\ell=min\{max\{k_1,...,k_{s+1-\ell}\},max\{k_{\ell},...,k_{s}\}\}$ for $\ell=1,2,...,s$.
\end{lem}

{\bf Proof:}
The proof is similar to that of Lemma \ref{lemNil1}.
Denote $P_j=\Pi_{j=1}^{\tau_1+...+\tau_j}A^{i_j}$.
Clearly, for any $A^{i_1},A^{i_2},...,A^{i_{\tau_1}}$,
\begin{eqnarray}
P_1&=\left(
  \begin{array}{cccccc}
    A_{11}^{\tau_1} & *  &* &*& \cdots  & *  \\
    \textbf{0} & A_{22}^{\tau_1} &* &*& \cdots & * \\
    \textbf{0} & \textbf{0} & A_{33}^{\tau_1} & *&\cdots &* \\
    \vdots & \vdots &\vdots&\ddots &\vdots&  \vdots \\
    \textbf{0} & \cdots & \textbf{0} & \textbf{0} &A_{s-1s-1}^{\tau_1}& *\\
    \textbf{0} & \cdots & \textbf{0} & \textbf{0}& \textbf{0} &A_{ss}^{\tau_1} \\
  \end{array}
\right)&\\\label{eqnBloNi2}
&=\left(
  \begin{array}{cccccc}
    \textbf{0} & *  &* &*&\cdots  & *  \\
    \textbf{0} & \textbf{0} &* &* &\cdots & * \\
    \textbf{0} & \textbf{0} &\textbf{0} &* & \cdots & * \\
    \vdots & \vdots &\vdots& \ddots & \vdots & \vdots\\
    \textbf{0} & \cdots & \textbf{0} & \textbf{0}& \textbf{0}&*\\
    \textbf{0} & \cdots & \textbf{0} & \textbf{0}& \textbf{0} & \textbf{0}\\
  \end{array}
\right)&\label{eqnBloNi3}
\end{eqnarray}

Denote the block matrix in the $(i,i+1)$ position of $P_1$ by
$B_i,i=1,...,s-1$. For any $i=1,...,s-1$ and any
$A^{i_{\tau_1+1}},A^{i_2},...,A^{i_{\tau_1+\tau_2}}$, if
$max\{k_2,...,k_{s}\}\}<max\{k_1,...,k_{s-1}\}$, then the block
matrix in the $(i,i+1)$ position of $P_2$ is
$B_iA_{(i+1)(i+1)}^{\tau_2}$, which equals to $\textbf{0}$; if
$max\{k_2,...,k_{s}\}\}\geq max\{k_1,...,k_{s-1}\}$, then the block
matrix in the $(i,i+1)$ position of $P_2$ is $A_{ii}^{\tau_2}B_i$,
which equals to $\textbf{0}$. Since the sequence is
arbitrary, the change of index is irrelevant for the case of
$max\{k_2,...,k_{s}\}\}\geq max\{k_1,...,k_{s-1}\}$  because we can
consider $A^{i_{\tau_2+1}}...A^{i_{\tau_1+\tau_2}}$ first, which
takes the form of (\ref{eqnBloNi3}), and then consider
$A^{i_1}A^{i_2}...A^{i_{\tau_1+\tau_2}}$. Besides, it is clear that
all the zero matrices in $P_1$ remain unchanged in $P_2$. Therefore,
for any $A^{i_1},...,A^{i_{\tau_1+\tau_2}}$, $P_2$ has the following
form:
\begin{eqnarray}
P_2&=\left(
  \begin{array}{ccccccc}
    \textbf{0} & \textbf{0} &* &* &*&\cdots  & *  \\
    \textbf{0} & \textbf{0} &\textbf{0} &*&* &\cdots & * \\
    \textbf{0} & \textbf{0} &\textbf{0} &\textbf{0} &*&\cdots & * \\
    \vdots & \vdots &\vdots& \ddots & \vdots & \vdots& \vdots\\
        \textbf{0} & \cdots & \textbf{0} & \textbf{0}& \textbf{0}&\textbf{0}&*\\
    \textbf{0} & \cdots & \textbf{0} & \textbf{0}& \textbf{0}&\textbf{0}&\textbf{0}\\
    \textbf{0} & \cdots & \textbf{0} & \textbf{0}& \textbf{0} & \textbf{0}&\textbf{0}\\
  \end{array}
\right)&\label{eqnBloNi4}
\end{eqnarray}
Following the same argument, we can prove that as $j$ increases,
$P_j$ has fewer nonzero block matrices in the upper-right positions.
Finally, for any sequence $A^{i_1},A^{i_2},...,A^{i_T}$ where
$T=\sum_{j=1}^s\tau_j$,
we have
$$
A^{i_1}A^{i_2}...A^{i_T}=\textbf{0}_{ns\times ns}
$$
The lemma is thus proved.
\hfill$\Box$

\begin{cor}\label{corNil}
Suppose that a matrix $A$ over $\F$ has the form of (\ref{eqnBloNil}), where
$A_{j}(j=1,...,s)$ are nilpotent matrices with nilpotent degree $n$.
Then $A$ is nilpotent with degree upper bounded by $ns$.
\end{cor}

\subsection{Theorems}

The following theorem is the main result of this section.

\begin{thm}\label{thmCon}
Suppose that $\bar \G$ is DAG and $A$ is not nilpotent. Then system
(\ref{eqnFollow}) achieve consensus with (\ref{eqnLeader}) using
control (\ref{eqnControl}) if and only if (i) $(A,b)$ is
stabilizable; (ii) there is $d\in\{1,...,p-1\}$ such that $d_i\equiv
d$ for any $i=1,...,N$.
\end{thm}

{\bf Proof:} Recall equation (\ref{eqnCan}) and let $K=(k_1,...,k_m)$, $K^c=KQ^{-1}=(k'_1,...,k'_m)$.

\emph{Necessity:} If (\ref{eqnFollow})  achieve
consensus with (\ref{eqnLeader}) under control (\ref{eqnControl}), then
$A^c-d_i b^cK^c$ and $A^{uc}$ are all nilpotent matrices for $i=1,...,N$, which means
that condition (i) holds. 
Condition (ii) will be proved by contradiction. 
Since $A^c-d_i b^cK^c$ is nilpotent, $det(\lambda I_s-(A^c-d_i
b^cK^c))=\lambda^s+(d_ik'_{s}-a_s)\lambda^{s-1}+...+(d_ik'_{2}-a_2)\lambda+
(d_ik'_{1}-a_1)=\lambda^s$. That is, $d_ik'_{l}- a_l=0$
for $i=1,...,N,l=1,...,s$. If there exists $\ell \in\{1,...,N\}$ such that $d_{\ell}=0$, then $a_l=0$ for $l=1,...s$. Then $A^c$ is
itself nilpotent and $A$ is therefore nilpotent, which  contradicts
with the assumption. If there exist $i_1,i_2\in\{1,...,N\}$ such that $d_{i_1}\neq d_{i_2}$, then $(d_{i_1}-d_{i_2})k'_l=0$ for any $l=1,...,s$. Because $\F$ has no zero divisor that is not $0$, $k'_l=0$ for $l=1,...,s$. This also implies that $a_l=0$ for $l=1,...s$, which contradicts with the assumption.
Thus, condition (ii) holds and the necessity part is proved.

\emph{Sufficiency:} Under conditions given in the theorem, we can
find a constant matrix $K$ such that $A-d_i bK$ is nilpotent for any
$i=1,...,N$. By Lemma \ref{lemNil1} and Corollary \ref{corNil}, the matrix $I_N\otimes
A+(\bar \A-\bar \D)\otimes bK$ is nilpotent and (\ref{connil})
holds. Then system (\ref{eqnFollow}) achieve
consensus with (\ref{eqnLeader}).\hfill$\Box$

\begin{rem}
The necessity proof of Theorem \ref{thmCon} implies that $(A^c, d_i b^c)$ have to be simultaneously stabilizable by a matrix $K$ for
$i=1,...,N$. To achieve that,
$d_i$ must be nonzero and equal to each other.
Moreover, the leader must be the
neighbor of each follower representing the source node of $\bar\G$,
which means the leader must be globally reachable in $\G$.
\end{rem}

\begin{rem}
To make (\ref{eqnFollow}) achieve consensus with (\ref{eqnLeader}),
$\bar \G$  is not necessary to be DAG. Generally speaking, matrices
$K$ and $\bar\A$ may be designed by solving a set of multi-variable
polynomial equations over $\F$, which was proved to be NP-hard
\cite{fra79}. 
\end{rem}



Finally, let us consider that the interaction graphs are time-varying. Let
$\{\G_p:p\in \mathcal{P}\}$ be the set of possible directed graphs
on node $\{0,1,...,N\}$, and $\{\bar\G_p=(\V_p,\E_p):p\in
\mathcal{P}\}$ be the set of induced subgraphs on node
$\{1,...,N\}$, where $\mathcal{P}=\{1,...,q\}$. The dependence of
the graphs upon time is determined by a discrete-time switching
signal $\sigma:\mathbb{Z}_+\rightarrow \mathcal{P}$, and the
underlying graph at instance $k$ is $\G_{\sigma(k)}
(\bar\G_{\sigma(k)})$. Let $d_i^k$ be the in-degree of agent $i$
under $\bar\G_{\sigma(k)}$.


\begin{thm}\label{thmSwi}
Suppose that $A$ is not nilpotent and the union of subgraphs
$\cup_{k\geq 1}\bar\G_{\sigma(k)}$ is DAG. If $(A,b)$ is
stabilizable and there is $d\in\{1,...,p-1\}$ such that $d_i^k\equiv
d$ for any $i=1,...,N$ and any $k$, then, under arbitrary switching signals, system (\ref{eqnFollow})
achieve consensus with (\ref{eqnLeader}) using control
(\ref{eqnControl}).
\end{thm}

{\bf Proof:} Find a constant matrix $K$ such that $A-d_i bK$ is
nilpotent for any $i=1,...,N$. If $\cup_{k\geq 1}\bar\G_{\sigma(k)}$
is DAG, then matrices $I_N\otimes A+(\bar \A_{\sigma(k)}-\bar
\D_{\sigma(k)})\otimes bK$ can be simultaneously transformed into
upper triangle forms by the same permutation matrix $P$. The conclusion can be proved easily by Lemma
\ref{lemNil2}. \hfill$\Box$

\begin{rem}
In both static and time-varying graph cases, the design for matrix
$K$ requires knowledge of the interaction topology and the system
pair $(A,b)$. But after that, each agent only needs to know the
relative information (between itself and its neighbors) and the edge weights (with its neighbors). In this
sense, (\ref{eqnControl}) can be considered as ``distributed''.
\end{rem}

\begin{rem}
Suppose that the leader has an external input
\begin{equation*}\label{eqnLead2}
  x_0(t+1) = Ax_0(t)+\bar bv(t)
\end{equation*}
Then the leader is able to present some pre-specified dynamical patterns
in $\F$ via static state feedback $v(t)=\bar Kx_0(t)$\cite{reg04}.
Provided that (\ref{eqnFollow}) achieve consensus with
(\ref{eqnLeader}), all the agents will exhibit the same dynamics no
matter what the initial conditions are. This can be seen as a method
to achieve quantized consensus if the states of agents are coded
somehow in $\F$.
\end{rem}

\subsection{Example}
Consider finite field $\mathbb F_3$ and equation (\ref{eqnFollow}) with
$$
A=\left(\begin{array}{ccccc}
    0&0&1&1&1\\
    2&0&0&1&2\\
    0&2&2&2&0\\
    0&0&1&1&2\\
    2&0&1&2&2
  \end{array}\right),
b=\left(\begin{array}{c}
    1 \\
    1 \\
    2 \\
    2 \\
    1
  \end{array}\right).
$$
There is an invertible matrix $Q$ such that $QAQ^{-1},Qb$ are in control companion forms where
$$
Q=\left(
    \begin{array}{ccccc}
     0&1&0&1&0\\
    2&0&1&2&1\\
    2&2&1&2&2\\
    2&2&2&0&2\\
    1&2&0&1&1
    \end{array}
  \right),
    $$
$$
  QAQ^{-1}=
  \left(
    \begin{array}{ccccc}
      0 & 1 & 0 & 0 & 0 \\
      0 & 0 & 1 & 0 & 0 \\
      0 & 0 & 0 & 1 & 0 \\
      2 & 1 & 0 & 2 & 2 \\
      0 & 0 & 0 & 0 & 0 \\
    \end{array}
  \right),
  Qb=\left(
  \begin{array}{c}
    0 \\
    0 \\
    0 \\
    1 \\
    0
  \end{array}\right).
$$
Since $P_A(\lambda)=\lambda(\lambda^4+2\lambda^3+\lambda+2)$ with
$\lambda^4+2\lambda^3+\lambda+2$ irreducible over $\mathbb F_3$, $A$
is not nilpotent. It can be calculated that system (\ref{eqnLeader}) has 4 cycles with length
$1,20,20,20$ \cite{tol05}. It is also easy to check that $(A,b)$ is
stabilizable.

Suppose there are $5$ agents and the possible weighted interaction
graphs $\G_1,\G_2$ are shown in Fig.1, with the weights
shown over edges. Because $\bar\G_1\cup\bar\G_2$ is DAG and the
in-degree for each agent is $1$ (noting that $2+2\equiv 1(mod\;
3)$), conditions of Theorem \ref{thmSwi} are satisfied. Choose
$K=[2,1,2,0,1]$ such that $A-bK$ is nilpotent. Then under control
protocol (\ref{eqnControl}) and arbitrary switching signals
$\sigma(k)\in\{1,2\}$, the system will achieve consensus after
finite time. Define $e_i=\sum_{j=1}^5|x_i^j-x_0^j|,i=1,...,4,$ as the error
between $x_i$ and $x_0$ where the arithmetics used are standard.
Given arbitrary initial conditions, typical evolutions of $e_i$ are
shown in Fig.2. We can find that $e_i$ approach 0
after a finite time, which indicates that (\ref{eqnFollow}) achieve
consensus with (\ref{eqnLeader}).

\begin{figure}\label{fig1}
  \centering
  \includegraphics[width=.7\linewidth]{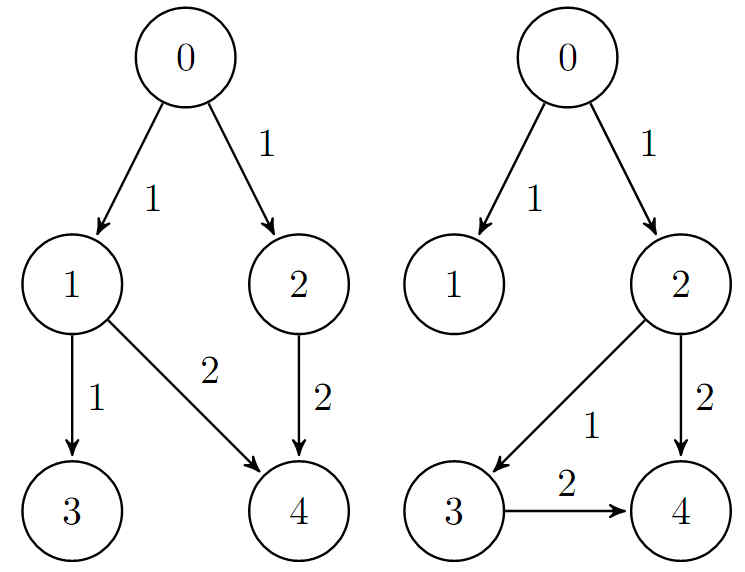}\\
 \caption{Interaction graphs $\G_1$ and $\G_2$}
\end{figure}

\begin{figure}\label{fig2}
  \centering
  \includegraphics[width=.7\linewidth]{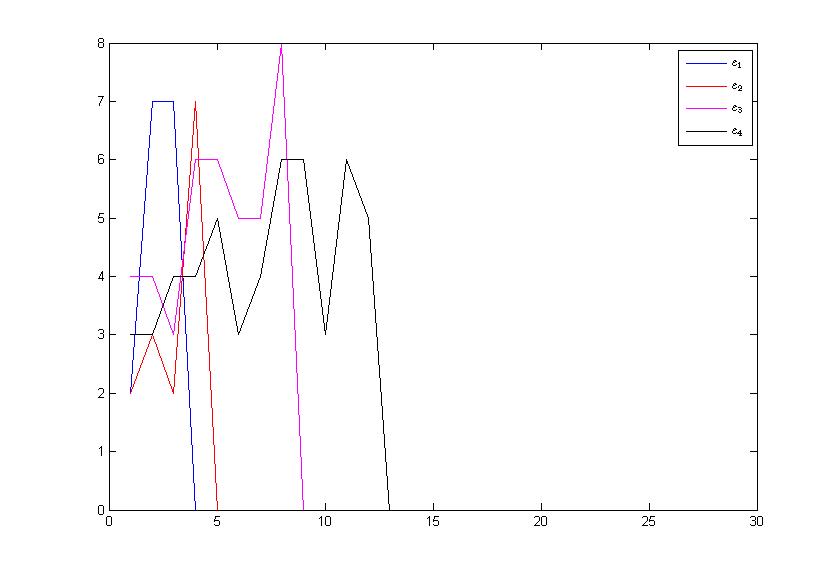}\\
 \caption{Errors between follower $i$ and the leader}
\end{figure}

\section{Conclusions}
In this paper, we formulated a leader-following consensus problem of
multi-agent systems over finite fields.  Then we gave sufficient
and/or necessary conditions for the agents to achieve the
leader-following consensus.  More general cases including general
graph topology and multiple inputs for the consensus in finite
fields are under investigation.



\end{document}